\documentclass[12pt,a4paper]{article}

\usepackage[latin2]{inputenc}
\usepackage{amssymb}
\usepackage{paralist}
\usepackage{amsmath, calc}
\usepackage{amsfonts}
\usepackage{amsthm}
\usepackage{fullpage}
\usepackage{enumerate}
\newtheorem{thm}{Theorem}

\newtheorem{prob}{Problem}

\newtheorem{lem}{Lemma}

\begin{document}

\title{Generalizations of some results about the regularity properties of an additive representation function}
\author{S\'andor Z. Kiss \thanks{Institute of Mathematics, Budapest
University of Technology and Economics, H-1529 B.O. Box, Hungary;
kisspest@cs.elte.hu;
This author was supported by the National Research, Development and Innovation Office NKFIH Grant No. K115288.}, Csaba
S\'andor \thanks{Institute of Mathematics, Budapest University of
Technology and Economics, H-1529 B.O. Box, Hungary, csandor@math.bme.hu.
This author was supported by the OTKA Grant No. K109789. This paper was supported
by the J\'anos Bolyai Research Scholarship of the Hungarian Academy of Sciences.} 
}
\date{}
\maketitle

\begin{abstract}
\noindent Let $A = \{a_{1},a_{2},\dots{}\}$  $(a_{1} < a_{2} < \dots{})$
    be an infinite sequence of nonnegative  integers, and let $R_{A,2}(n)$ denote
    the number of solutions of $a_{x}+a_{y}=n$  $(a_{x},a_{y}\in A)$. 
P. Erdős, A. Sárközy and V. T. Sós proved  that if
$\lim_{N\to\infty}\frac{B(A,N)}{\sqrt{N}}=+\infty$ then
$|\Delta_{1}(R_{A,2}(n))|$ cannot be bounded, where $B(A,N)$ denotes
the number of blocks formed by consecutive integers in $A$ up to $N$
and $\Delta_{l}$ denotes the $l$-th difference. Their result was extended 
to $\Delta_{l}(R_{A,2}(n))$ for any fixed $l\ge2$. In this paper we give further generalizations of this problem.


\end{abstract}
\section{Introduction}

Let $\mathbb{N}$ denote the set of nonnegative integers. Let $k \ge 2$ be a fixed integer and let $A = \{a_{1}, a_{2}, \dots{}\}$  $(a_{1} < a_{2} < \dots{})$ be an infinite sequence of nonnegative
integers. For $n = 0, 1, 2, \dots{}$ let $R_{A,k}(n)$ denote the number of
solutions of $a_{i_{1}} + a_{i_{2}} + \dots + a_{i_{k}} = n$,  $a_{i_{1}} \in A, \dots ,a_{i_{k}} \in A$, and we put
\[
A(n) = \sum_{\overset{a \in A}{a \le n}}1.
\]
\noindent We denote the cardinality of a set $H$ by $\#H$. Let $B(A,N)$ denote the number of blocks formed by consecutive integers in $A$ up to $N$, i.e.,
\[
B(A, N) = \sum_{\overset{n \le N}{n \in A, n-1 \notin A}}1.
\]
If $s_{0}, s_{1}, \dots{}$ is
given sequence of real numbers then let $\Delta_{l}s_{n}$ denote the $l$-th
difference of the sequence $s_{0}, s_{1}, s_{2}, \dots{}$ defined by 
$\Delta_{1}s_{n} = s_{n+1}-s_{n}$ and $\Delta_{l}s_{n} 
= \Delta_{1}(\Delta_{l-1}s_{n})$.

\noindent In a series of papers \cite{ES}, \cite{EK}, \cite{ER} P. Erdős, 
A. Sárközy and 
V.T. Sós studied the regularity properties of the function
$R_{A,2}(n)$. In \cite{ER} they proved the following theorem:

\smallskip

\noindent{\bf  Theorem A}~{\em
If $\lim_{N\to\infty}\frac{B(A,N)}{\sqrt{N}} = \infty$, then $|\Delta_{1}(R_{A,2}(n))| = |R_{A,2}(n+1)-R_{A,2}(n)|$ cannot be bounded.}

\smallskip

\noindent In \cite{ER} they also showed that the above result is nearly best 
possible:

\smallskip

\noindent{\bf  Theorem B}~{\em
For all $\varepsilon > 0$, there exists an infinite sequence $A$
such that
\begin{itemize}
\item[(i)] $B(A,N) \gg N^{1/2-\varepsilon}$,
\item[(ii)] $R_{A,2}(n)$ is bounded so that also $\Delta_{1}R_{A,2}(n)$ is 
bounded. 
\end{itemize}
}
\noindent Recently, \cite{RK} A. S\'ark\"ozy extended the above results in the finite set of residue classes modulo a fixed $m$. 

In \cite{SI} Theorem A was extended to any  $k > 2$ :

\smallskip

\noindent{\bf  Theorem C}~{\em
If $k \ge 2$ is an integer and $\lim_{N\to\infty}\frac{B(A,N)}{\sqrt[k]{N}} 
= \infty$, and  $l \le k$, then $|\Delta_{l}R_{A,k}(n)|$ cannot be bounded.
}

\smallskip

\noindent It was shown \cite{SQ} that the above result is nearly best
possible. 

\smallskip

\noindent{\bf  Theorem D}~{\em
For all $\varepsilon > 0$, there exists an infinite sequence $A$
such that
\begin{itemize}
\item[(i)] $B(A,N) \gg N^{1/k-\varepsilon}$,
\item[(ii)] $R_{A,k}(n)$ is bounded so that also $\Delta_{l}R_{A,k}(n)$ is 
bounded if $l \le k$. 
\end{itemize}
}

\noindent In this paper we consider $R_{A,2}(n)$, thus simply write $R_{A,2}(n) = R_{A}(n)$. A set of positive integers $A$ is called Sidon set if $R_{A}(n) \le 2$. Let $\chi_{A}$ denote the characteristic function of the set $A$, i.e.,
\[
\chi_{A}(n) =
\left\{
\begin{aligned}
1 \textnormal{, if } n \in A \\
0 \textnormal{, if } n \notin A.
\end{aligned} \hspace*{3mm}
\right.
\]
Let $\lambda_{0}, \dots{}, \lambda_{d}$ be arbitrary integers with $\sum_{i=0}^{d}|\lambda_{i}| > 0$ Let $\underline{\lambda} = (\lambda_{0}, \dots{}, \lambda_{d})$ and define the function 
\[
B(A, \underline{\lambda}, n) = \left|\left \{m: m \le n, \sum_{i=0}^{d}\lambda_{i}\chi_{A}(m - i) \ne 0\right \}\right|.
\]

\begin{thm}
We have
\[
\limsup_{n \rightarrow \infty}\left|\sum_{i=0}^{d}\lambda_{i}R_{A}(n - i)\right| \ge
\limsup_{n \rightarrow \infty}\frac{|\sum_{i=0}^{d}\lambda_{i}|}{2(d+1)^{2}}\left(\frac{B(A, \underline{\lambda}, n)}{\sqrt{n}}\right)^{2}.
\]
\end{thm}
\noindent The next theorem shows that the above result is nearly best possible: 
\begin{thm}
Let $\sum_{i=0}^{d}\lambda_{i} > 0$.
Then for every positive integer $N$ there exists a set $A$ such that
\[
\limsup_{n \rightarrow \infty}\left|\sum_{i=0}^{d}\lambda_{i}R_{A}(n - i)\right| \le \limsup_{n \rightarrow \infty}4\sum_{i=0}^{d}|\lambda_{i}|\left(\frac{B(A, \underline{\lambda}, n)}{\sqrt{n}}\right)^{2}
\]
and
\[
\limsup_{n \rightarrow \infty}\frac{B(A, \underline{\lambda}, n)}{\sqrt{n}} \ge N.
\]
\end{thm}

\begin{thm}
Let $\sum_{i=0}^{d}\lambda_{i} = 0$. Then we have
\[
\limsup_{n \rightarrow \infty}\left|\sum_{i=0}^{d}\lambda_{i}R_{A}(n - i)\right| \ge \limsup_{n \rightarrow \infty}\frac{\sqrt{2}}{e^{2}\sum_{i=0}^{d}|\lambda_{i}|}
\frac{B(A, \underline{\lambda}, n)}{\sqrt{n}}.
\]
\end{thm}
It is easy to see that if $\underline{\lambda} = (\lambda_{0}, \lambda_{1}) = (-1, 1)$ then $B(A, \underline{\lambda}, n) \ge B(A,n)$ thus Theorem 3 implies 
Theorem A. It is natural to ask that the exponent of $\frac{B(A, \underline{\lambda}, n)}{\sqrt{n}}$ in the right hand side can be improved.
                                                                                                                                                             
\begin{prob}
Is it true that if $\sum_{i=0}^{d}\lambda_{i} = 0$ then there exists a positive
constant $C(\underline{\lambda})$ depends only on $\underline{\lambda}$ 
such that for every set of nonnegative integers $A$ we have
\[
\limsup_{n \rightarrow \infty}\left|\sum_{i=0}^{d}\lambda_{i}R_{A}(n - i)\right| \ge \limsup_{n \rightarrow \infty}C(\underline{\lambda})\cdot \left(\frac{B(A, \underline{\lambda}, n)}{\sqrt{n}}\right)^{3/2}?
\]
\end{prob}
\noindent In the next theorem we prove that the exponent cannot grow over $3/2$.
\begin{thm}
Let $\sum_{i=0}^{d}\lambda_{i} = 0$.
For every positive integer $N$ there exists a set $A \subset \mathbb{N}$ such 
that
\[
N \le \limsup_{n \rightarrow \infty}\left(\frac{B(A, \underline{\lambda}, n)}{\sqrt{n}}\right)
< \infty
\]
and
\[
\limsup_{n \rightarrow \infty}\left|\sum_{i=0}^{d}\lambda_{i}R_{A}(n - i)\right| \le \limsup_{n \rightarrow \infty}48(d+1)^42^{3d+7.5}\sum_{i=0}^d|\lambda _i|\left(\frac{B(A, \underline{\lambda}, n)}{\sqrt{n}}\right)^{3/2}\left(\log \frac{B(A, \underline{\lambda}, n)}{\sqrt{n}}\right)^{1/2}.
\]
\end{thm}

\section{Proof of Theorem 1}

Since $-\underline{\lambda} = (-\lambda_{0}, \dots{}, -\lambda_{d})$ and clearly 
\[
\limsup_{n \rightarrow \infty}\left|\sum_{i=0}^{d}\lambda_{i}R_{A}(n - i)\right| = \limsup_{n \rightarrow \infty}\left|\sum_{i=0}^{d}(-\lambda_{i})R_{A}(n - i)\right|,
\]
$B(A, \underline{\lambda}, n) = B(A, -\underline{\lambda}, n)$, therefore
\[
\limsup_{n \rightarrow \infty}\frac{|\sum_{i=0}^{d}\lambda_{i}|}{2(d+1)^{2}}\left(\frac{B(A, \underline{\lambda}, n)}{\sqrt{n}}\right)^{2} = \limsup_{n \rightarrow \infty}\frac{|\sum_{i=0}^{d}(-\lambda_{i})|}{2(d+1)^{2}}\left(\frac{B(A, -\underline{\lambda}, n)}{\sqrt{n}}\right)^{2},
\]
thus we may assume that $\sum_{i=0}^{d}\lambda_{i} > 0$. On the other hand we may suppose that
\[
\limsup_{n \rightarrow \infty}\left(\frac{B(A, \underline{\lambda}, n)}{\sqrt{n}}\right)^
{2} > 0.
\]
It follows from the definition of the limsup that there exists a sequence $n_{1}, n_{2}, \dots{}$ such that
\[
\lim_{j \rightarrow \infty}\frac{B(A, \underline{\lambda}, n_{j})}{\sqrt{n_{j}}} = \limsup_{n \rightarrow \infty}\frac{B(A, \underline{\lambda}, n)}{\sqrt{n}}.
\]

To prove Theorem 1 we give a lower and an upper estimation to
\begin{equation}\label{eq:1}
\sum_{\sqrt[3]{n_{j}} < n \le 2n_{j}}\left(\sum_{i=0}^{d}\lambda_{i}R_{A}(n-i)\right).
\end{equation}

The comparison of the two bounds will give the result.
First we give an upper estimation. Clearly we have
\[
\left |\sum_{\sqrt[3]{n_{j}} < n \le 2n_{j}}\left(\sum_{i=0}^{d}\lambda_{i}R_{A}(n-i)
\right)\right | \le
\sum_{\sqrt[3]{n_{j}} < n \le 2n_{j}}\left |\sum_{i=0}^{d}\lambda_{i}R_{A}(n-i) \right |
\]
\[
\le 2n_{j}\text{max}_{\sqrt[3]{n_{j}} < n \le 2n_{j}}\left |\sum_{i=0}^{d}\lambda_{i}R_{A}(n-i)\right |.
\]
In the next step we give a lower estimation to (\ref{eq:1}). It is clear that
\[
\sum_{\sqrt[3]{n_{j}} < n \le 2n_{j}}\left(\sum_{i=0}^{d}\lambda_{i}R_{A}(n-i)\right)
= \sum_{\sqrt[3]{n_{j}} < n \le 2n_{j}}(\lambda_{0} + \dots{} + \lambda_{d})R_{A}(n)
\]
\[
- \Big((\lambda_{1} + \dots{} + \lambda_{d})R_{A}(2n_{j}) + (\lambda_{2} + \dots{} + \lambda_{d})R_{A}(2n_{j}-1) + \lambda_{d}R_{A}(2n_{j}-d+1)\Big)
\]
\[
+(\lambda_{1} + \dots{} + \lambda_{d}) R_{A}(\lfloor \sqrt[3]{n_{j}} \rfloor) + (\lambda_{2} + \dots{} + \lambda_{d}) R_{A}(\lfloor \sqrt[3]{n_{j}} \rfloor - 1) + \dots{} + \lambda_{d}R_{A}(\lfloor \sqrt[3]{n_{j}} \rfloor - d + 1).
\]
Obviously,
\[
R_{A}(m) = \#\{(a, a^{'}): a + a^{'} = m, a, a^{'} \in A\} \le 2\cdot \#\{(a, a^{'}): a + a^{'} = m,  a \le a^{'}, a, a^{'} \in A\}
\]
\[
\le 2\cdot \#\{(a: a \le m/2, a \in A\} = 2A(m/2).
\]
It follows that
\[
\sum_{\sqrt[3]{n_{j}} < n \le 2n_{j}}\left(\sum_{i=0}^{d}\lambda_{i}R_{A}(n-i)\right) \ge (\lambda_{0} + \dots{} + \lambda_{d})\sum_{\sqrt[3]{n_{j}} < n \le 2n_{j}}R_{A}(n) - \left(\sum_{i=0}^{d}|\lambda_{i}|\right)2A(n_{j})2d
\]
\[
\ge \left(\sum_{i=0}^{d}\lambda_{i}\right)\#\{(a,a^{'}): a+a^{'} = n, \sqrt[3]{n_{j}} < a,a^{'} \le n_{j}, a,a^{'} \in A \} - \left(\sum_{i=0}^{d}|\lambda_{i}|\right)4dA(n_{j})
\]
\[
= \left(\sum_{i=0}^{d}\lambda_{i}\right)(A(n_{j}) - A(\sqrt[3]{n_{j}}))^{2} - O(A(n_{j})).
\]
The inequaltity $\sum_{i=0}^{d}\lambda_{i}\chi_{A}(m - i) \ne 0$ implies that $[m-d, m] \cap A \ne 0$. Then we have $\{m: m \le n, \sum_{i=0}^{d}\lambda_{i}\chi_{A}(m - i) \ne 0\} \subseteq
\cup_{a \le n, \\ a \in A}[a, a + d]$, which implies that $B(A, \underline{\lambda}, n) \le \left |\cup_{a \le n, \\ a \in A}[a, a + d]\right | \le A(n)(d + 1)$.
By the definition of $n_{j}$ there exists a constant $c_{1}$ such that
\[
\frac{B(A, \underline{\lambda}, n_{j})}{\sqrt{n_{j}}} > c_{1} > 0.
\]
It follows that $A(n_{j}) > \frac{c_{1}}{d+1}\sqrt{n_{j}}$ and clearly $\sqrt[3]{n_{j}} \ge A(\sqrt[3]{n_{j}})$. By using these facts we get that
\[
\left(\sum_{i=0}^{d}\lambda_{i}\right)(A(n_{j}) - A(\sqrt[3]{n_{j}}))^{2} - O(A(n_{j})) = (1+o(1))\left(\sum_{i=0}^{d}\lambda_{i}\right)A(n_{j})^{2} \ge 
\]
\[
(1+o(1))\left(\sum_{i=0}^{d}\lambda_{i}\right)\frac{B(A, \underline{\lambda}, n_{j})^{2}}{(d+1)^{2}}.
\]
Comparing the lower and the upper estimations we get that
\[
2n_{i}\max_{\sqrt[3]{n_{j}} < n \le 2n_{j}}\left|\sum_{i=0}^{d}\lambda_{i}R_{A}(n-i)\right| \ge \sum_{\sqrt[3]{n_{j}} < n \le 2n_{j}}\left(\sum_{i=0}^{d}\lambda_{i}R_{A}(n-i)\right)
\]
\[
\ge (1+o(1))\frac{\sum_{i=0}^{d}\lambda_{i}}{(d+1)^{2}}B^{2}(A, \underline{\lambda}, n_{j}),
\]
this implies that
\begin{equation}\label{eq:2}
\max_{\sqrt[3]{n_{j}} < n \le 2n_{j}}\left|\sum_{i=0}^{d}\lambda_{i}R_{A}(n-i)\right| \ge (1+o(1))\frac{\sum_{i=0}^{d}\lambda_{i}}{2(d+1)^{2}}\left(\frac{B(A, \underline{\lambda}, n_{j})}{\sqrt{n_{j}}}\right)^{2}.
\end{equation}
To complete the proof we distinguish two cases.
When
\[
\limsup_{n \rightarrow \infty}\left(\frac{B(A, \underline{\lambda}, n)}{\sqrt{n}}\right)^{2} < \infty
\]
then
\[
\max_{\sqrt[3]{n_{j}} < n \le 2n_{j}}\left|\sum_{i=0}^{d}\lambda_{i}R_{A}(n-i)\right| \ge (1+o(1))\frac{\sum_{i=0}^{d}\lambda_{i}}{2(d+1)^{2}}\left(\frac{B(A, \underline{\lambda}, n_{j})}{\sqrt{n_{j}}}\right)^{2}
\]
\[
= (1+o(1)) \frac{\sum_{i=0}^{d}\lambda_{i}}{2(d+1)^{2}}\limsup_{n \rightarrow \infty}\left(\frac{B(A, \underline{\lambda}, n)}{\sqrt{n}}\right)^{2},
\]
which gives the result.

When
\[
\limsup_{n \rightarrow \infty}\left(\frac{B(A, \underline{\lambda}, n)}{\sqrt{n}}\right)^{2} = \infty
\]
then 
\[
\limsup_{j \rightarrow \infty}\left(\frac{B(A, \underline{\lambda}, n_{j})}{\sqrt{n_{j}}}\right)^{2} = \infty,
\]
which implies by (\ref{eq:2}) that $\limsup_{n \rightarrow \infty}\left|\sum_{i=0}^{d}\lambda_{i}R_{A}(n-i)\right| = \infty$, which gives the result.

\section{Proof of Theorem 2}

It is well known \cite{HR} that there exists a Sidon set $S$ with
\[
\limsup_{n \rightarrow \infty}\frac{S(n)}{\sqrt{n}} \ge \frac{1}{\sqrt{2}},
\]
where $S(n)$ denotes the number of elements of $S$ up to $n$. Define the set $T$ by removing the elements $s$ and $s^{'}$ from $S$
when $s - s^{'} \le (N + 1)(d + 1)$. It is clear that $T(n) \ge S(n) - 2(N+1)(d+1)$ and define the set $A$ by
\[
A = T \cup (T + (d + 1)) \cup (T + 2(d + 1)) \cup \dots{} \cup (T + N(d + 1)).
\]
It is easy to see that $A(n) \ge (N + 1)T(n) - N$. We will prove that
$B(A, \underline{\lambda}, n) \ge A(n) - d$. By the definitions of the sets $T$ and $A$ we get that if $a < a^{'}$, $a$, $a^{'}\in A$ then $a - a^{'} \ge d + 1$. If
\[
\sum_{i=0}^{d}\lambda_{i}\chi_{A}(m - i) \ne 0
\]
then there is exactly one term, which is nonzero. Fix an index $w$ such that
$\lambda_{w} \ne 0$. It follows that $\sum_{i=0}^{d}\lambda_{i}\chi_{A}(a+w-i) \ne 0$ for every $a\in A$. Hence,
\[
|B(A, \underline{\lambda}, n)| \ge \#\{a: a + w \le n, a\in A\} = A(n-w) \ge A(n) - w \ge A(n) - d
\]
\[
\ge (N + 1)T(n) - N - d \ge (N+1)S(n) - 2(N+1)^{2}(d + 1) - N - d.
\]
Thus we have
\[
\frac{B(A, \underline{\lambda}, n)}{\sqrt{n}} \ge (N + 1)\frac{S(n)}{\sqrt{n}} - \frac{2(N+1)^{2}(d + 1) + N + d}{\sqrt{n}}
\]
and
\[
\limsup_{n \rightarrow \infty}\left(\frac{B(A, \underline{\lambda}, n)}{\sqrt{n}}\right)^{2} \ge \frac{(N+1)^{2}}{2}\ge N.
\]
By the definition of $A$, we have
\[
R_{A}(m) = \sum_{i=0}^{N}\sum_{j=0}^{N}\#\{(t,t^{'}): (t+i(d+1)) + (t+j(d+1)) = m, t,t^{'} \in T\}
\]
\[
= \sum_{i=0}^{N}\sum_{j=0}^{N}R_{T}(m-(i+j)(d+1)) \le 2(N+1)^{2}.
\]
Then we have
\[
\left|\sum_{i=0}^{d}\lambda_{i}R_{A}(n - i)\right| \le \left(\sum_{i=0}^{d}|\lambda_{i}|\right)
\max_{n}R_{A}(n) \le 2\left(\sum_{i=0}^{d}|\lambda_{i}|\right)(N+1)^{2}\le 
\]
\[
\le \limsup_{n \rightarrow \infty}4\cdot \left(\sum_{i=0}^{d}|\lambda_{i}|\right)\left(\frac{B(A, \underline{\lambda}, n)}{\sqrt{n}}\right)^{2},
\]
which gives the result.

\section{Proof of Theorem 3}

In the first case we assume that
\[
\limsup_{n \rightarrow \infty}\frac{\sqrt{2}}{e^{2}\sum_{i=0}^{d}|\lambda_{i}|}
\frac{B(A, \underline{\lambda}, n)}{\sqrt{n}} < \infty.
\]

We prove by contradiction. Assume that contrary to the conclusion of Theorem 3
we have

\begin{equation}\label{eq:3}
\limsup_{n \rightarrow \infty}\left|\sum_{i=0}^{d}\lambda_{i}R_{A}(n - i)\right| < \limsup_{n \rightarrow \infty}\frac{\sqrt{2}}{e^{2}\sum_{i=0}^{d}|\lambda_{i}|}\frac{B(A, \underline{\lambda}, n)}{\sqrt{n}}.
\end{equation}

Throughout the remaining part of the proof of Theorem 3 we use the following
notations: $N$ denotes a positive integer. We write $e^{2i\pi\alpha} = e(\alpha)$
and we put
$r = e^{-1/N}$, $z = re(\alpha)$ where $\alpha$ is a real variable (so that
a function of form $p(z)$ is a function of the real variable
${\alpha: p(z) = p(re(\alpha)) = P(\alpha)}$). We write
$f(z) = \sum_{a \in A}z^{a}$. (By $r < 1$, this infinite series
and all the other infinite series in the remaining part of the proof are
absolutely convergent).\\

\indent We start out from the integral
$I(N) = \int\limits_0^1|f(z)(\sum_{i=0}^{d}\lambda_{i}z^{i})|^{2}d\alpha$. We will give lower and upper bound
for $I(N)$. The comparison of these bounds will give a contradiction.

First we will give a lower bound for $I(N)$. We write
\[
f(z)\left(\sum_{i=0}^{d}\lambda_{i}z^{i}\right) = \left(\sum_{n=0}^{\infty}\chi_{A}(n)z^{n}\right)\left(\sum_{i=0}^{d}\lambda_{i}z^{i}\right)
\]
\[
= \sum_{n=0}^{\infty}(\lambda_{0}\chi_{A}(n) + \lambda_{1}\chi_{A}(n-1) + \dots{} + \lambda_{d}\chi_{A}(n-d))z^{n}.
\]
It is clear that if $\lambda_{0}\chi_{A}(n) + \lambda_{1}\chi_{A}(n-1) + \dots{} + \lambda_{d}\chi_{A}(n-d) \ne 0$, then $(\lambda_{0}\chi_{A}(n) + \lambda_{1}\chi_{A}(n-1) + \dots{} + \lambda_{d}\chi_{A}(n-d))^{2} \ge 1$. Thus, by the Parseval formula, we have
\[
I(N) = \int_{0}^{1}\left|f(z)\left(\sum_{i=0}^{d}\lambda_{i}z^{i}\right)\right|^{2}
d\alpha
\]
\[
= \int_{0}^{1}\left|\sum_{n=0}^{\infty}(\lambda_{0}\chi_{A}(n) + \lambda_{1}\chi_{A}(n-1) + \dots{} + \lambda_{d}\chi_{A}(n-d))z^{n}\right|^{2}d\alpha
\]
\[
= \sum_{n=0}^{\infty}\left(\lambda_{0}\chi_{A}(n) + \lambda_{1}\chi_{A}(n-1) + \dots{} + \lambda_{d}\chi_{A}(n-d)\right)^{2}r^{2n} \ge e^{-2}\sum_{\overset{n \le N}{\lambda_{0}\chi_{A}(n) + \lambda_{1}\chi_{A}(n-1) + \dots{} + \lambda_{d}\chi_{A}(n-d) \ne 0}}1
\]
\[
= e^{-2}B(A, \underline{\lambda}, N).
\]
Now we will give an upper bound for $I(N)$. Since
the sums $\sum_{i=0}^{d}|\lambda_{i}R_{A}(n-i)|$ are nonnegative integers it follows from
(\ref{eq:3}) that there exists an $n_{0}$ and an $\varepsilon > 0$ such that
\begin{equation}\label{eq:4}
\sum_{i=0}^{d}\left|\lambda_{i}R_{A}(n-i)\right| \le \limsup_{n \rightarrow \infty}\frac{\sqrt{2}}{e^{2}\sum_{i=0}^{d}|\lambda_{i}|}\frac{B(A, \underline{\lambda}, n)}{\sqrt{n}}(1 - \varepsilon).
\end{equation}
for every $n > n_{0}$. On the other hand there exists an infinite sequence of real numbers
$n_{0} < n_{1} < n_{2} < \dots{} < n_{j} < \dots{}$ such that
\[
\limsup_{n \rightarrow \infty}\frac{\sqrt{2}}{e^{2}\sum_{i=0}^{d}|\lambda_{i}|}\frac{B(A, \underline{\lambda}, n)}{\sqrt{n}}\sqrt{1 - \varepsilon} < \frac{\sqrt{2}}{e^{2}\sum_{i=0}^{d}|\lambda_{i}|}\frac{B(A, \underline{\lambda}, n_{j})}{\sqrt{n_{j}}}.
\]
We get that
\begin{equation}\label{eq:5}
\limsup_{n \rightarrow \infty}\frac{\sqrt{2}}{e^{2}\sum_{i=0}^{d}|\lambda_{i}|}\frac{B(A, \underline{\lambda}, n)}{\sqrt{n}}(1 - \varepsilon) < \frac{\sqrt{2}}{e^{2}\sum_{i=0}^{d}|\lambda_{i}|}\frac{B(A, \underline{\lambda}, n_{j})}{\sqrt{n_{j}}}\sqrt{1 - \varepsilon}.
\end{equation}
Obviously, $f^{2}(z) = \sum_{n = 0}^{\infty}R_{A}(n)z^{n}$.
By our indirect
assumption, the Cauchy inequality and the Parseval formula we have
\[
I(N) = \int\limits_0^1\left|f(z)\left(\sum_{i=0}^{d}\lambda_{i}z^{i}\right)\right|^{2}d\alpha
\le \left(\sum_{i=0}^{d}|\lambda_{i}|\right)\int\limits_0^1\left|f^{2}(z)\left(\sum_{i=0}^{d}\lambda_{i}z^{i}\right)\right|d\alpha
\]
\[
= \left(\sum_{i=0}^{d}|\lambda_{i}|\right)\int\limits_0^1\left|\left(\sum_{n = 0}^{\infty}R_{A}(n)z^{n}\right)\left(\sum_{i=0}^{d}\lambda_{i}z^{i}\right)\right|d\alpha = \left(\sum_{i=0}^{d}|\lambda_{i}|\right)\int\limits_0^1\left|\sum_{n = 0}^{\infty}\left(\sum_{i=0}^{d}\lambda_{i}R_{A}(n-i)\right)z^{n}\right|d\alpha
\]
\[
\le \left(\sum_{i=0}^{d}|\lambda_{i}|\right) \left (\int\limits_0^1\left|\sum_{n = 0}^{\infty}\left(\sum_{i=0}^{d}\lambda_{i}R_{A}(n-i)\right)z^{n}\right|^{2}d\alpha \right )^{1/2} = \left(\sum_{i=0}^{d}|\lambda_{i}|\right) \left (\sum_{n = 0}^{\infty}(\sum_{i=0}^{d}\lambda_{i}R_{A}(n-i))^{2}r^{2n} \right )^{1/2}.
\]
In view of (\ref{eq:4}), (\ref{eq:5}) and the lower bound for $I(n_{j})$ we 
\[
e^{-2}B(A, \underline{\lambda}, n_{j}) < I(n_{j}) < \left(\sum_{i=0}^{d}|\lambda_{i}|\right) \left (\sum_{n = 0}^{\infty}\left(\sum_{i=0}^{d}\lambda_{i}R_{A}(n-i)\right)^{2}r^{2n} \right )^{1/2}
\]
\[
\le \left(\sum_{i=0}^{d}|\lambda_{i}|\right)\left (\sum_{n = 0}^{n_{0}}\left(\sum_{i=0}^{d}\lambda_{i}R_{A}(n-i)\right)^{2}r^{2n} + \sum_{n = n_{0}+1}^{\infty}\left(\frac{\sqrt{2}}{e^{2}\sum_{i=0}^{d}|\lambda_{i}|}\frac{B(A, \underline{\lambda}, n_{j})}{\sqrt{n_{j}}}\sqrt{1 - \varepsilon}  \right)^{2}r^{2n}\right )^{1/2}
\]
\[
< \left(\sum_{i=0}^{d}|\lambda_{i}|\right)\left (c_{2} + \sum_{n = 0}^{\infty}\left(\frac{2}{e^{4}(\sum_{i=0}^{d}|\lambda_{i}|)^{2}}\frac{B^{2}(A, \underline{\lambda}, n_{j})}{n_{j}}(1 - \varepsilon)  \right)r^{2n}\right )^{1/2},
\]
where $c_{2}$ is a constant. Taking the square of both sides we get that
\begin{equation}\label{eq:7}
e^{-4}B^{2}(A, \underline{\lambda}, n_{j}) < \left(\sum_{i=0}^{d}|\lambda_{i}|\right)^{2}\left (c_{2} + \frac{2}{e^{4}(\sum_{i=0}^{d}|\lambda_{i}|)^{2}}\frac{B^{2}(A, \underline{\lambda}, n_{j})}{n_{j}}(1 - \varepsilon)\sum_{n=0}^{\infty}r^{2n}\right).
\end{equation}

It is easy to see that
\[
1 - e^{-x} = x - \frac{x^2}{2!} + \frac{x^3}{3!} - \dots > x - \frac{x^2}{2!} =
x(1-\frac{x}{2}) > \frac{x}{x+1}
\]
for $0 < x < 1$. Applying this observation, where $r = e^{-1/n_{j}}$ we have
\[
\sum_{n=0}^{\infty}r^{2n} = \frac{1}{1 - r^{2}} = \frac{1}{1-e^{-\frac{2}{n_{j}}}}
\]
\[
< \frac{n_{j}}{2} + 1.
\]
In view of (\ref{eq:7}) we obtain that
\[
e^{-4}B^{2}(A, \underline{\lambda}, n_{j}) < \left(\sum_{i=0}^{d}|\lambda_{i}|\right)^{2}\left (c_{2} + \frac{2}{e^{4}(\sum_{i=0}^{d}|\lambda_{i}|)^{2}}\frac{B^{2}(A, \underline{\lambda}, n_{j})}{n_{j}}(1 - \varepsilon)\left(\frac{n_{j}}{2} + 1 \right) \right)
\]
\[
< c_{3} + e^{-4}B^{2}(A, \underline{\lambda}, n_{j})(1 - \varepsilon),
\]
where $c_{3}$ is an absolute constant and it follows that
\[
B^{2}(A, \underline{\lambda}, n_{j}) < c_{3}e^{4} + B^{2}(A, \underline{\lambda}, n_{j})(1 - \varepsilon),
\]
or in other words
\[
B^{2}(A, \underline{\lambda}, n_{j}) < \frac{c_{3}e^{4}}{\varepsilon},
\]
which is a contradiction if $n_{j}$ is large enough because $\lim_{j \rightarrow \infty}B(A, \underline{\lambda}, n_{j}) = \infty$. This proves the first case.

\noindent Assume that
\[
\limsup_{n \rightarrow \infty}\frac{\sqrt{2}}{e^{2}\sum_{i=0}^{d}|\lambda_{i}|}
\frac{B(A, \underline{\lambda}, n)}{\sqrt{n}} = \infty.
\]
Then there exists a sequence $n_{1} < n_{2} < \dots{}$ such that
\[
\limsup_{j \rightarrow \infty}\frac{B(A, \underline{\lambda}, n_{j})}{\sqrt{n_{j}}} = \infty.
\]
We prove by contradiction. Suppose that
\[
\limsup_{n \rightarrow \infty}\left|\sum_{i=0}^{d}\lambda_{i}R_{A}(n - i)\right| < \infty.
\]
Then there exists a positive constant $c_{4}$ such that $|\sum_{i=0}^{d}\lambda_{i}R_{A}(n - i)| < c_{4}$ for every $n$.
It follows that
\[
e^{-2}B(A, \underline{\lambda}, n_{j}) < I(n_{j}) < \left(\sum_{i=0}^{d}|\lambda_{i}|\right) \left (\sum_{n = 0}^{\infty}\left(\sum_{i=0}^{d}\lambda_{i}R_{A}(n-i)\right)^{2}r^{2n} \right )^{1/2} < \left(c_{4}\sum_{n = 0}^{\infty}r^{2n}\right )^{1/2}
< c_{5}\sqrt{n_{j}},
\]
thus we have
\[
\frac{B(A, \underline{\lambda}, n_{j})}{\sqrt{n_{j}}} < c_{5}e^{2},
\]
where $c_{5}$ is a positive constant, which is absurd.

\section{Proof of Theorem 4}

We argue as S\'ark\"ozy in \cite{RK}.
In the first step we will prove the following lemma:

\begin{lem}
There exists a set $C_{M} \subset [0, M(d+1) - 1]$ for which $|R_{C_{M}}(n) - R_{C_{M}}(n-1)| \le
12 \sqrt{M(d+1)\log M(d+1)}$ for every nonnegative integer $n$ and $B(C_{M}, \underline{\lambda}, M(d+1)-1) \ge \frac{M}{2^{d+2}}$ if $M$
is large enough.
\end{lem}

\noindent {\bf Proof of Lemma 1} To prove the lemma we use the probabilistic method due to Erd\H{o}s and R\'enyi. There is an excellent summary about this method in books \cite{AS} and \cite{HR}. Let $\mathbb{P}(E)$ denote the probability of an event $E$ in a probability space and let $\mathbb{E}(X)$ denote the expectation of a random variable $X$.
Let us define a random set $C$ with
$\mathbb{P}(n \in C) = \frac{1}{2}$ for every $0 \le n \le M(d + 1) - 1$. In the first step we show that
\[
\mathbb{P}\left(\text{max}_{n}|R_{C}(n) - R_{C}(n-1)| > 12 \sqrt{M(d+1)\log M(d+1)}\right) < \frac{1}{2}.
\]
Define the indicator random variable
\[
\varrho_{C}(n) =
\left\{
\begin{aligned}
1 \textnormal{, if } n \in C \\
0 \textnormal{, if } n \notin C.
\end{aligned} \hspace*{3mm}
\right.
\]
It is clear that
\[
R_{C}(n) = 2\sum_{k<n/2}\varrho_{C}(k)\varrho_{C}(n-k) + \varrho_{C}(n/2)
\]
sum of independent indicator random variables. Define the random variable
$\zeta_{i}$ by

$\zeta_{i} = \varrho_{C}(i)\varrho_{C}(n-i)$. Then we have
\[
R_{C}(n) = 2X_{n} + Y_{n},
\]
where $X_{n} = \zeta_{0} + \dots{} + \zeta_{\lfloor \frac{n-1}{2} \rfloor}$ and $Y_{n} = \varrho_{C}(n/2)$.

\noindent \textbf{Case 1.} Assume that $0 \le n \le M(d+1) - 1$. Obviously, $\mathbb{P}(\zeta_{i} = 0) = \frac{3}{4}$ and $\mathbb{P}(\zeta_{i} = 1) = \frac{1}{4}$ and
\[
\mathbb{E}(X_{n}) = \frac{\lfloor \frac{n+1}{2} \rfloor}{4}.
\]
As $Y_{n} \le 1$, it is easy to see that the following events satisfy the following relations
\begin{eqnarray*}
&&\{\text{max}_{0 \le n \le M(d+1)-1}|R_{C}(n) - R_{C}(n-1)| > 12\sqrt{M(d+1)\log M(d+1)}\}\\
&\subseteq& \left \{\text{max}_{0 \le n \le M(d+1)-1}\left|R_{C}(n) - \frac{n}{4} + R_{C}(n-1) -  \frac{n-1}{4}\right| > 10\sqrt{M(d+1)\log M(d+1)}\right \}\\
&\subseteq& \left \{\text{max}_{0 \le n \le M(d+1)-1}\left(\left|R_{C}(n) - \frac{n}{4}\right| + \left|R_{C}(n-1) -  \frac{n-1}{4}\right|\right) > 10\sqrt{M(d+1)\log M(d+1)}\right \}\\
&\subseteq& \left\{\text{max}_{0 \le n \le M(d+1)-1}\left|R_{C}(n) - \frac{n}{4}\right|  > 5\sqrt{M(d+1)\log M(d+1)}\right\}\\
&=& \left \{\text{max}_{0 \le n \le M(d+1)-1}\left|2X_{n} + Y_{n} - \frac{n}{4}\right| > 5\sqrt{M(d+1)\log M(d+1)}\right \}\\
&\subseteq& \left \{\text{max}_{0 \le n \le M(d+1)-1}\left|2X_{n} - \frac{n}{4}\right| > 4\sqrt{M(d+1)\log M(d+1)}\right\}\\
&=& \left \{\text{max}_{0 \le n \le M(d+1)-1}\left|X_{n} - \frac{n}{8}\right| > 2\sqrt{M(d+1)\log M(d+1)}\right \}\\
&\subseteq& \left\{\text{max}_{0 \le n \le M(d+1)-1}\left|X_{n} - \frac{\lfloor \frac{n+1}{2} \rfloor}{4}\right| > \sqrt{M(d+1)\log M(d+1)}\right\}.
\end{eqnarray*}
It follows that
\[
\mathbb{P}\left(\text{max}_{0 \le n \le M(d+1)-1}|R_{C}(n) - R_{C}(n-1)| > 12\sqrt{M(d+1)\log M(d+1)}\right) 
\]
\[
\le \mathbb{P}\left(\text{max}_{0 \le n \le M(d+1)-1}\left|X_{n} - \frac{\lfloor \frac{n+1}{2} \rfloor}{4}\right| > \sqrt{M(d+1)\log M(d+1)}\right).
\]
\[
\le \sum _{n=0}^{M(d+1)-1} \mathbb{P}\left(\left|X_{n} - \frac{\lfloor \frac{n+1}{2} \rfloor}{4}\right| > \sqrt{M(d+1)\log M(d+1)}\right).
\]
It follows from the Chernoff type bound \cite{AS}, Corollary A 1.7. that if the random variable $X$ has Binomial distribution with parameters $m$ and $p$ then for $a > 0$ we have
\begin{equation}\label{eq:6}
\mathbb{P}(|X - mp| > a) \le 2e^{-2a^{2}/m}.
\end{equation}
Applying (\ref{eq:6}) to  $\lfloor \frac{n+1}{2} \rfloor$ and $p = \frac{1}{4}$ we have
\begin{equation}\label{eq:100}
\mathbb{P}\left(\left|X_{n} - \frac{\lfloor \frac{n+1}{2} \rfloor}{4}\right| > \sqrt{M(d+1)\log M(d+1)}\right) < 2\cdot \text{exp}\left(\frac{-2M(d+1)\log M(d+1)}{\lfloor \frac{n+1}{2} \rfloor}\right) 
\end{equation}
\[
\le 2e^{-4\frac{M(d+1)\log M(d+1)}{M(d+1)}} = 2e^{-4\log M(d+1)} = \frac{2}{(M(d+1))^{4}} < \frac{1}{4M(d+1)}.
\]
It follows that
\begin{equation}\label{eq:8}
 \mathbb{P}(\{\text{max}_{0 \le n \le M(d+1)-1}|R_{C}(n) - R_{C}(n-1)| > 12\sqrt{M(d+1)\log M(d+1)}\})
\end{equation}
\[
< \frac{M(d+1)}{4M(d+1)} = \frac{1}{4}.
\]

\noindent \textbf{Case 2.} Assume that $M(d+1) \le n \le 2M(d+1) - 2$.

Obviously, $\mathbb{P}(\zeta_{i} = 0) = \frac{3}{4}$ and
$\mathbb{P}(\zeta_{i} = 1) = \frac{1}{4}$ when $n - M(d+1) < i < \frac{n}{2}$, and if $0 \le i \le n - M(d+1)$ then $\zeta_{i} = 0$. Clearly we have
\[
\mathbb{E}(X_{n}) = \frac{\lfloor \frac{2M(d+1)-1-n}{2} \rfloor}{4}.
\]
As $Y_{n} \le 1$, it is easy to see that the following relations holds among the events
\begin{align*}
&\{\text{max}_{M(d+1) \le n \le 2M(d+1)-2}|R_{C}(n) - R_{C}(n-1)| > 12\sqrt{M(d+1)\log M(d+1)}\}\\
\subseteq& \Bigg \{\max_{M(d+1)\le n \le 2M(d+1)-2}\left|R_{C}(n) - \frac{M(d+1)-\frac{n}{2}}{4} + R_{C}(n-1) -  \frac{M(d+1)-\frac{n-1}{2}}{4}\right|\\
>& 10\sqrt{M(d+1)\log M(d+1)}\Bigg \}\\
\subseteq& \Bigg \{\max_{M(d+1) \le n \le 2M(d+1)-2}\left(\left|R_{C}(n) - \frac{M(d+1)-\frac{n}{2}}{4}\right| + \left|R_{C}(n-1) - \frac{M(d+1)-\frac{n-1}{2}}{4}\right|\right)\\ 
>& 10\sqrt{M(d+1)\log M(d+1)}\Bigg \}\\
\subseteq& \left \{\text{max}_{M(d+1)-1 \le n \le 2M(d+1)-2}\left|R_{C}(n) - \frac{M(d+1)-\frac{n}{2}}{4}\right| > 5\sqrt{M(d+1)\log M(d+1)}\right \}\\
=& \left \{\text{max}_{M(d+1)-1 \le n \le 2M(d+1)-2}\left|2X_{n} + Y_{n} - \frac{2M(d+1)-n}{4}\right| > 5\sqrt{M(d+1)\log M(d+1)}\right \}\\
\subseteq& \left\{\text{max}_{M(d+1)-1 \le n \le 2M(d+1)-2}\left|2X_{n} - \frac{2M(d+1)-n}{4}\right| > 4\sqrt{M(d+1)\log M(d+1)}\right \}\\
=& \left \{\text{max}_{M(d+1)-1 \le n \le 2M(d+1)-2}\left|X_{n} - \frac{2M(d+1)-n}{8}\right| > 2\sqrt{M(d+1)\log M(d+1)}\right \}\\
\subseteq& \left \{\text{max}_{M(d+1)-1 \le n \le 2M(d+1)-2}\left|X_{n} - \frac{\lfloor \frac{2M(d+1)-1-n}{2} \rfloor}{4}\right| > \sqrt{M(d+1)\log M(d+1)}\right \}\\
\end{align*}
It follows that
\[
\mathbb{P}\left(\text{max}_{M(d+1)-1 \le n \le 2M(d+1)-2}|R_{C}(n) - R_{C}(n-1)| > 12\sqrt{M(d+1)\log M(d+1)}\right)
\]
\[
\le \mathbb{P}\left(\max _{M(d+1)-1\le n \le 2M(d+1)-1}\left|X_{n} - \frac{\lfloor \frac{2M(d+1)-1-n}{2} \rfloor}{4}\right| > \sqrt{M(d+1)\log M(d+1)}\right)
\]
\[
\le \sum _{n=M(d+1)-1}^{2M(d+1)-2}\mathbb{P}\left(\left|X_{n} - \frac{\lfloor \frac{2M(d+1)-1-n}{2} \rfloor}{4}\right| > \sqrt{M(d+1)\log M(d+1)}\right)
\]

Applying (\ref{eq:6}) for $m = \frac{\lfloor \frac{2M(d+1)-1-n}{2} \rfloor}{4}$
and $p = \frac{1}{4}$ we have for $M(d+1) \le n \le 2M(d+1)-2$
\[
\mathbb{P}\left(\left|X_{n} - \frac{\lfloor \frac{2M(d+1)-1-n}{2} \rfloor}{4}\right| > \sqrt{M(d+1)\log M(d+1)}\right) < 2\cdot \text{exp}\left(\frac{-2M(d+1)\log M(d+1)}{\lfloor \frac{2M(d+1)-1-n}{2} \rfloor}\right) 
\]
\[
< 2e^{-4\frac{M(d+1)\log M(d+1)}{M(d+1)}} = 2e^{-4\log M(d+1)} 
= \frac{2}{(M(d+1))^{4}}< \frac{1}{4M(d+1)}
\]
and by (\ref{eq:100}) we have
\[
\mathbb{P}\left(\left|X_{M(d+1)-1} - \frac{\lfloor \frac{n+1}{2} \rfloor}{4}\right| > \sqrt{M(d+1)\log M(d+1)}\right) < \frac{1}{4M(d+1)}.
\]
It follows that
\begin{equation}\label{eq:200}
\mathbb{P}\left(\text{max}_{M(d+1) \le n \le 2M(d+1)-2}|R_{C}(n) - R_{C}(n-1)| > 12\sqrt{M(d+1)\log M(d+1)}\right)\end{equation} 
\[
< \frac{M(d+1)}{4M(d+1)} = \frac{1}{4}.
\]

By (\ref{eq:8}) and (\ref{eq:200}) we get that

\begin{equation}\label{eq:800}
\mathbb{P}\left(\text{max}_{0 \le n \le 2M(d+1)-2}|R_{C}(n) - R_{C}(n-1)| > 12\sqrt{M(d+1)\log M(d+1)}\right) < \frac{1}{2}.
\end{equation}

In the next step we show that
\[
\mathbb{P}\left(B(C, \underline{\lambda}, M(d+1)-1) < \frac{M}{2^{d+2}}\right) < \frac{1}{2}.
\]
It is clear that the following events $E_{1}, \dots{}, E_{M}$ are independent:
\begin{align*}
E_{1} & = \left \{\sum_{i=0}^{d}\lambda_{i}\varrho_{C}(d-i) \ne 0\right \},\\
E_{2} & = \left \{\sum_{i=0}^{d}\lambda_{i}\varrho_{C}(d+1+d-i) \ne 0\right \},\\
&\mathrel{\makebox[\widthof{=}]{\vdots}}\\
E_{M} & = \left \{\sum_{i=0}^{d}\lambda_{i}\varrho_{C}((m-1)(d+1)+d-i) \ne 0\right 
\}.
\end{align*}
Obviously, $\mathbb{P}(E_{i}) = \mathbb{P}(E_{j})$, where $1 \le i,j \le M$.
Let $p = \mathbb{P}(E_{1})$. It is clear that there exists an index $u$ such that $\lambda_{u} \ne 0$. Thus we have
\[
p \ge \mathbb{P}(\varrho_{C}(0) = 0, \varrho_{C}(1) = 0, \dots{}, \varrho_{C}(u-1) = 0, \varrho_{C}(u) = 1, \varrho_{C}(u+1) = 0, \dots{}, \varrho_{C}(d) = 0)
\]
\[
= \frac{1}{2^{d+1}}.
\]
Define the random variable $Z$ as the number of occurrence of the events $E_{j}$.
It is easy to see that $Z$ has Binomial distribution with parameters $M$ and $p$. Apply the Chernoff bound (\ref{eq:6}) we get that
\[
\mathbb{P}\left(|Z - Mp| > \frac{Mp}{2}\right) < 2e^{\frac{-2(Mp/2)^{2}}{M}} < 2e^{-\frac{M}{2}\cdot 2^{-2d-2}}
< \frac{1}{2}
\]
if $M$ is large enough. On the other hand, we have
\[
\frac{1}{2} > \mathbb{P}\left(|Z - Mp| > \frac{Mp}{2}\right) 
\ge \mathbb{P}\left(Z < \frac{Mp}{2}\right) \ge \mathbb{P}\left(Z < \frac{M}{2^{d+2}}\right).
\]
Hence,
\begin{equation}\label{eq:9}
\mathbb{P}\left(B(C,\underline{\lambda}, 2M(d+1)-2) < \frac{M}{2^{d+2}}\right) < \frac{1}{2}.
\end{equation}
Let $\mathcal{E}$ and $\mathcal{F}$ be the events
\[
\mathcal{E} = \left \{\max_{0 \le n \le 2M(d+1)-2}|R_{C}(n) - R_{C}(n-1)| > 12\sqrt{M(d+1)\log M(d+1)}\right \},
\]
\[
\mathcal{F} = \left \{B(C,\underline{\lambda}, M(d+1)-1) < \frac{M}{2^{d+2}}\right \}.
\]
It follows from (\ref{eq:800}) and (\ref{eq:9}) that
\[
\mathbb{P}\left(\mathcal{E} \cup \mathcal{F}\right) < 1,
\]
then
\[
\mathbb{P}\left(\overline{\mathcal{E}} \cap \overline{\mathcal{F}} \right) > 0,
\]
therefore there exists a suitable set $C_{M}$ if $M$ is large enough, which completes the proof of Lemma 1.

We are ready to prove Theorem 4.
It is well known \cite{HR} that there exists a Sidon set $S$ with
\[
\limsup_{n \rightarrow \infty}\frac{S(n)}{\sqrt{n}} \ge \frac{1}{\sqrt{2}},
\]
where $S(n)$ is the number of elements of $S$ up to $n$. Let $s, s^{'} \in S$ and assume that $s > s^{'}$. Define $S_{M} = S \setminus \{s, s^{'} \in S: s - s^{'} \le 2M(d+1)\}$ and let
$A = C_{M} + S_{M}$, where $C_{M}$ is the set from the lemma.
\[
\left|\sum_{i=0}^{d}\lambda_{i}R_{A}(n - i)\right| = \left|\sum_{i=0}^{d}\lambda_{i}\#\{(a,a^{'}): a + a^{'} = n-i, a,a^{'} \in A\}\right|
\]
\[
= \left|\sum_{i=0}^{d}\lambda_{i}\#\{(s,c,s^{'},c^{'}): s + c + s^{'} + c^{'} = n-i, s,s^{'} \in S_{M}, c,c^{'} \in C_{M}\}\right|
\]
\[
= \left|\sum_{i=0}^{d}\sum_{j=0}^{2M(d+1)}\lambda_{i}\#\{(s,c,s^{'},c^{'}): c + c^{'} = j,
s + s^{'} = n - i - j, s,s^{'} \in S_{M}, c,c^{'} \in C_{M}\}\right|
\]
\[
= \left|\sum_{i=0}^{d}\sum_{j=0}^{2M(d+1)}\lambda_{i}R_{C_{M}}(j)R_{S_{M}}(n-i-j)\right|
\]
\[
= \left|\sum_{j=0}^{2M(d+1)}\sum_{i=0}^{d}\lambda_{i}R_{C_{M}}(j)R_{S_{M}}(n-i-j)\right| = \left|\sum_{k=0}^{2M(d+1)+d}\sum_{i=0}^{d}\lambda_{i}R_{C_{M}}(k-i)R_{S_{M}}(n-k)
\right|
\]
\[
\left|\sum_{k=0}^{2M(d+1)+d}R_{S_{M}}(n-k)\sum_{i=0}^{d}\lambda_{i}R_{C_{M}}(k-i)\right| \le
\sum_{k=0}^{2M(d+1)+d}R_{S_{M}}(n-k)\left|\sum_{i=0}^{d}\lambda_{i}R_{C_{M}}(k-i)\right|
\]
\[
\le 2(M+1)(d+1)2\cdot \max_{k}\left|\sum_{i=0}^{d}\lambda_{i}R_{C_{M}}(k-i)\right|.
\]
In the next step we give an upper estimation to $|\sum_{i=0}^{d}\lambda_{i}R_{C_{M}}(k-i)|$.
We have
\begin{align*}
&|\lambda_{0}R_{C_{M}}(k) + \dots{} + \lambda_{d}R_{C_{M}}(k-d)|\\
&= |\lambda_{0}(R_{C_{M}}(k)-R_{C_{M}}(k-1)) + (\lambda_{0} +\lambda_{1})(R_{C_{M}}(k-1)-R_{C_{M}}(k-2)) + \dots{} \\
&+ (\lambda_{0} +\lambda_{1} + \dots{} + \lambda_{d-1})(R_{C_{M}}(k-d+1)-R_{C_{M}}(k-d)) + (\lambda_{0} +\lambda_{1} + \dots{} + \lambda_{d})R_{C_{M}}(k-d)|.
\end{align*}
Since $\sum_{i=0}^{d}\lambda_{i} = 0$, the last term in the previous sum is zero. Then we have
\[
|\lambda_{0}R_{C_{M}}(k) + \dots{} + \lambda_{d}R_{C_{M}}(k-d)| \le d\left(\sum_{i=0}^{d}|\lambda_{i}|\right)\max_{t}|R_{C_{M}}(t)-R_{C_{M}}(t-1)|\le 
\]
\[
12d\sum_{i=0}^d|\lambda _i|\sqrt{M(d+1)\log M(d+1)}.
\]
Then we have
\[
\left|\sum_{i=0}^{d}\lambda_{i}R_{A}(n - i)\right| \le 48d\sum_{i=0}^d|\lambda _i|(M(d+1))^{3/2}(\log M(d+1))^{1/2}.
\]
We give a lower estimation to
\[
\limsup_{n \rightarrow \infty}\frac{B(A, \underline{\lambda}, n)}{\sqrt{n}}.
\]
If $0 \le v \le M(d+1) - 1$ and $\sum_{i=0}^{d}\lambda_{i}\chi_{C_{M}}(v - i) \ne 0$
then $\sum_{i=0}^{d}\lambda_{i}\chi_{A}(s+v-i) \ne 0$ for every $s \in S_{M}$. Then
we have
\[
B(A, \underline{\lambda}, n) \ge (S_{M}(N) - 1)B(C_{M}, \underline{\lambda}, (M+1)-1).
\]
Thus we have
\[
\limsup_{n \rightarrow \infty}\frac{B(A, \underline{\lambda}, n)}{\sqrt{n}} \ge \frac{M}{2^{d+2.5}}.
\]
It follows that
\[
\limsup_{n \rightarrow \infty}\left|\sum_{i=0}^{d}\lambda_{i}R_{A}(n - i)\right|
\le 48d\sum_{i=0}^d|\lambda _i|\left((M(d+1))^{3}\log M(d+1)\right)^{1/2}
\]
\[
\le \limsup_{n \rightarrow \infty}48(d+1)^42^{3d+7.5}\sum_{i=0}^d|\lambda _i|\left(\left(\frac{B(A, \underline{\lambda}, n)}{\sqrt{n}}\right)^{3}\log \frac{B(A, \underline{\lambda}, n)}{\sqrt{n}}\right)^{1/2},
\]
if $M$ is large enough. The proof of Theorem 4 is completed.

\end{document}